\documentclass{article}
\usepackage[utf8]{inputenc}
\usepackage{amsmath,amsfonts,amssymb,amsthm,mathtools} 
\usepackage{dsfont}
\usepackage{array}
\usepackage{tikz}
\usepackage{wrapfig}
\usepackage[normalem]{ulem}

\newtheorem{Lemma}[]{Lemma}
\newtheorem{Remark}[]{Remark}
\newtheorem{Theorem}[]{Theorem}

\def\Xint#1{\mathchoice
   {\XXint\displaystyle\textstyle{#1}}%
   {\XXint\textstyle\scriptstyle{#1}}%
   {\XXint\scriptstyle\scriptscriptstyle{#1}}%
   {\XXint\scriptscriptstyle\scriptscriptstyle{#1}}%
   \!\int}
\def\XXint#1#2#3{{\setbox0=\hbox{$#1{#2#3}{\int}$}
     \vcenter{\hbox{$#2#3$}}\kern-.5\wd0}}

\def\dashint{\Xint-}

\def\u{{\mathbf{u}}}

\title{Scale invariant bounds for mixing in the Rayleigh-Taylor instability}
\author{Konstantin Kalinin\footnote{Max Planck Institute for Mathematics in the Sciences, Inselstr. 22, 04103 Leipzig, Germany. E-mail: konstantin.kalinin@mis.mpg.de}, Govind Menon\footnote{Division of Applied Mathematics, Box F, Brown University, Providence, RI 02912, USA. E-mail: govind\_menon@brown.edu}, and Bian Wu\footnote{Max Planck Institute for Mathematics in the Sciences, Inselstr. 22, 04103 Leipzig, Germany. E-mail: bian.wu@mis.mpg.de}}
\date{}

\begin{document}
    \maketitle
    \begin{abstract}
    We study the Rayleigh-Taylor instability for two miscible, incompressible, inviscid fluids. 
    Scale-invariant estimates for the size of the mixing zone and coarsening of internal structures in the fully nonlinear regime are established following techniques introduced for the Saffman-Taylor instability in~\cite{OttoMenon2004}. These bounds provide optimal scaling laws and reveal the strong role of dissipation in slowing down mixing. 
    \end{abstract}

    \subsection*{Keywords.} Rayleigh-Taylor instability, Saffman-Taylor instability, fluid mixing. 
    \section{Introduction}

We consider an infinite column of two fluids of densities $\rho_+>\rho_->0$ in a gravitational field, initially arranged so that the heavy fluid is on top. The interface between the fluids is unstable and any small perturbation leads to turbulent mixing. This stability problem was first studied by Rayleigh~\cite{Rayleigh1883} and then analyzed again in Taylor's work~\cite{Taylor1950} (see also~\cite{Fermi-vonNeumann}). The Rayleigh-Taylor instability occurs in several natural settings (e.g. collapse of massive stars, formation of clouds) and technological applications (e.g. laser and electromagnetic implosions) (see~\cite{Sharp1984,Ye1,Ye2} for reviews). 

The following three stages of evolution are observed in experiments and numerical simulations. In an early stage, the growth of the instability is determined by linear analysis~\cite{Taylor1950}. This is followed by a fully nonlinear scaling regime where the turbulent mixing zone grows at the rate $\alpha_\pm Agt^2$, where $A = (\rho_+-\rho_-)/(\rho_++\rho_-)$ is the Atwood number, $g$ is the gravitational constant and $\alpha_\pm$ are empirically determined constants corresponding to the upper and lower edges of the mixing zone. This stage is also characterized by the formation of internal structures that have been termed bubbles and spikes. Finally, there is a late stage dominated by dissipation. In miscible and inviscid fluids, this dissipation is caused by molecular diffusion.

This paper provides rigorous results in the fully nonlinear regime and the late stage of the Rayleigh-Taylor instability that is based on the underlying partial differential equations. We use bulk quantities, energy balance, and an interpolation inequality to estimate the size of the mixing zone. The underlying techniques were introduced by one of the authors (GM) and Otto in~\cite{OttoMenon2004} for the study of gravity driven mixing for flows in a porous medium (the Saffman-Taylor instability). In this paper, we show that the underlying techniques may be adapted to the Rayleigh-Taylor instability, though several modifications of the main ideas are necessary. These include the inclusion of inertia in the energy balance, as well as a modification of our arguments to predict the coarsening rate of bubbles and spikes. We also include a modification of the scale-invariant energy balance to {understand the} crossover from the regime of turbulent mixing to the dissipative regime, when molecular diffusion acts to slow down mixing. 

Our results are consistent with a recent rigorous analysis by Gebhard {\em et al\/}~\cite{Gebhard2021}. Unlike this work, we do not construct weak solutions; rather we establish a priori estimates for classical solutions that must be satisfied by any weak solutions, following the approach in~\cite{OttoMenon2004}.  While there is an extensive literature on the early stage of the Rayleigh-Taylor instability~(see for example~\cite{Guo-Spirn,Guo}), to the best of our knowledge, the article~\cite{Gebhard2021} and this paper are the first rigorous results on the fully nonlinear regime that accord with experiments.

\subsection{The main model}
We model the evolution of the Rayleigh-Taylor instability for miscible, inviscid, incompressible fluids with the following system of PDEs
\begin{align}
    &\rho_t + \u \cdot \nabla \rho - \Delta \rho = 0, \label{01-transport}\\
    &\rho (\u_t + (\u \cdot \nabla) \u) - (\nabla \rho \cdot \nabla) \u + \nabla p + \rho g \mathbf{e}_z = 0, &\nabla \cdot \u = 0 \label{01-momentum}.
\end{align}
Here $\rho$ is the density, $\u$ is the velocity field, $g$ is the acceleration, $p$ is the pressure, and $\mathbf{e}_z$ is the unit vector in the vertical direction. Equations~\eqref{01-transport}-\eqref{01-momentum} are posed in an infinite strip $[0, L]^{d - 1} \times \mathbb{R}$, with $d\geq 2$. We denote points $x = (y, z)$ with $y \in [0,L]^{d-1}$ and $z \in \mathbb{R}$. The parameter $L$ plays the role of the Peclet number and quantifies the effect of miscibility. We assume periodic boundary conditions in the transverse direction $y$. 

The model consists of an advection-diffusion equation~\eqref{01-transport} for the density $\rho$ coupled with the incompressible Euler equations~\eqref{01-momentum} for the velocity $\u$ and pressure $p$. The diffusion term in equation~\eqref{01-transport} models miscibility. The inclusion of diffusion in equation~\eqref{01-transport} also requires the inclusion of an additional term in the momentum equation~\eqref{01-momentum}, corresponding to the transport driven by diffusion. This ensures energy balance as a fundamental law, as shown in Lemma~\ref{03-lm-energy-conservation}. Diffusion regularizes the evolution and causes the coarsening of bubbles and spikes as quantified in Theorem~\ref{02-thm-crossover}. 

We do not consider the effect of viscosity in this paper. The interplay between molecular diffusion and viscosity in the Rayleigh-Taylor instability is quite subtle, and several competing models have been proposed such as the variable density model (VDM). The VDM is based on Fickian diffusion which leads to artificial compressibility and, hence, to additional dissipation mechanisms (see~\cite{Cook-Dimotakis} for a review and~\cite{Gibbon} for an analysis). 

The steady, perfectly stratified state is
\begin{equation}
\label{eq:ic1}
    \u_0 = 0, \; \rho_0 = \begin{cases}
        \rho_+, &z \geqslant 0, \\
        \rho_-, &z < 0.
    \end{cases}
\end{equation}
This state serves as a reference configuration. Throughout this paper, the initial condition is a small perturbation of this state.

\begin{figure}
    \centering
    \begin{tikzpicture}[scale = 0.3]
    \draw  plot[smooth, tension=.7] coordinates {(-2.5,0) (-1,-0.5) (0.5,0) (1,1.5) (0.5,3.5) (-0.5,3) (-1.5,2.5) (-2.5,4) (-1.5,6.5) (0.5,8.5) (3,7.5) (4.5,6) (4,3.5) (3.5,2) (2,3.5) (2.5,1.5) (3.5,0) (5,-1) (5,-3.5) (5,-5) (4.5,-6) (3.5,-5.5) (3,-4) (2.5,-2.5) (2,-3) (2,-4) (2.5,-5.5) (3.5,-7) (5.5,-8) (7.5,-7.5) (8,-6) (8,-3.5) (7,-3) (7,-4.5) (6.5,-6) (6,-4.5) (5.5,-3) (6,-1) (7.5,0.5) (9,0.5)};
    \node (v1) at (-4.5,8.5) {};
    \node (v2) at (0.5,8.5) {};
    \draw  (v1) edge (v2);
    \node (v4) at (0.5,0) {};
    \node (v3) at (-4.5,0) {};
    \draw  (v3) edge (v4);
    \node (v5) at (-3.5,8.5) {};
    \node (v6) at (-3.5,0) {};
    \draw  (v5) edge (v6);
    \node (v8) at (-2.5,4) {};
    \node (v7) at (-2.5,10.5) {};
    \node (v9) at (4.5,6) {};
    \node (v10) at (4.5,10.5) {};
    \draw  (v7) edge (v8);
    \draw  (v9) edge (v10);
    \node (v11) at (-2.5,10) {};
    \node (v12) at (4.5,10) {};
    \draw  (v11) edge (v12);
    \node (v14) at (5.5,-8) {};
    \node (v13) at (-4.5,-8) {};
    \draw  (v13) edge (v14);
    \node (v15) at (-3.5,-8) {};
    \draw  (v6) edge (v15);
    \node at (2,-9.5) {};
    \node (v17) at (2,-9.5) {};
    \node (v19) at (8,-9.5) {};
    \node (v16) at (2,-4) {};
    \node (v18) at (8,-6) {};
    \draw  (v16) edge (v17);
    \draw  (v18) edge (v19);
    \node (v20) at (2,-9) {};
    \node (v21) at (8,-9) {};
    \draw  (v20) edge (v21);
    
    \node at (-5,4) {\rotatebox{90}{$a_+(t)$}};
    \node at (1,11) {$b_+(t)$};
    \node at (-5,-5) {\rotatebox{90}{$a_-(t)$}};
    \node at (5.5,-10) {$b_-(t)$};
\end{tikzpicture}
    \caption{Caricature of the main scales}
    \label{fig:scales}
\end{figure}
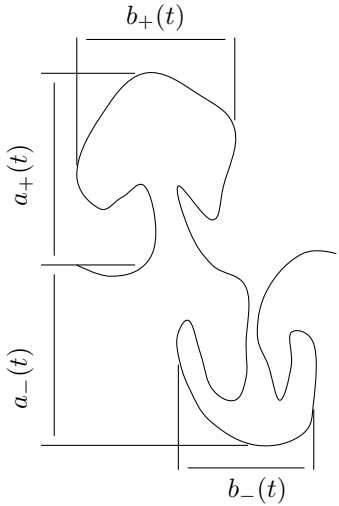

Let us outline our results informally. 
We define mixing lengths in the vertical and transverse directions using the potential energy and a mixing entropy and establish $L$-independent estimates for these quantities with the optimal scaling in time. Our main result is Theorem~\ref{02-main-thm}. This theorem provides $L$-independent estimates which may be interpreted as the one-sided bounds
\begin{equation}
    \label{eq:scale3}
     \alpha_- Agt^2 \leq a_-(t), \quad  a_+(t) \leq \alpha_{+} Agt^2, \quad b_{\pm}(t) \geq \beta_\pm \sqrt{t},
\end{equation}
where $a_{\pm}(t)$ measures the upper and lower edge of the mixing zone, $b_{\pm}(t)$ measures coarsening of structures within the mixing zone, and $\alpha_\pm$ and $\beta_\pm$ are universal constants depending on $A$ alone. See Figure~\ref{fig:scales} for a caricature of~$a_\pm$ and~$b_\pm$.

The proof of Theorem~\ref{02-main-thm} is based on two main arguments: energy balance (Section~\ref{sec:diff:ineq}) and general interpolation (Theorem~\ref{03-thm-interp} in Section~\ref{sec:interp:ineq}). In Theorem~\ref{03-thm-interp}, we bound the mixing entropy of profiles in terms of the potential energy in a scale-invariant manner. We also determine the sharp constant for this inequality and show that the optimal profile corresponds to the rarefaction wave solution for a scalar conservation law. This profile may be interpreted as the coarse-grained evolution of a mixing zone. It is the comparison with scalar conservation laws that allows us to determine explicit formulas for the constants $\alpha_\pm$. These constants are compared with experimental and numerical results in the area (see Section~\ref{sec:rarefactions}) in order to illustrate the phenomenon of dissipative slowdown.


The methods used to establish these bounds follow earlier work~\cite{OttoMenon2004} on the Saffman-Taylor instability.Here we have the following system of PDEs
\begin{align}
    & s_t + \u \cdot \nabla s - \Delta s = 0, \label{01-st-transport}\\
    &\u + \nabla p = - s \mathbf{e}_z, & \nabla \cdot \u = 0 \label{01-st-darcy},
\end{align}
in the same geometry, where $s \in [0, 1]$ is saturation. The initial data is given by a small perturbation of the perfectly stratified interface
\begin{equation}\label{01-st-stratified-profile}
    s_0(z) = \begin{cases}
        1, & z \geqslant 0, \\
        0, & z < 0.
    \end{cases}
\end{equation}
The system~\eqref{01-st-transport}--\eqref{01-st-darcy} with a small perturbation of~\eqref{01-st-stratified-profile} as initial datum describes the evolution of the instability arising in a porous medium reservoir or in a Hele-Shaw cell under the influence of gravity.

The phenomenon of dissipative slowdown for the Saffman-Taylor instability was analyzed in~\cite{OttoMenon2004,OttoMenon2006} as follows. The use of $L$-independent energy balance as in this work shows that the mixing layer grows as $2t$~\cite[Theorem 1 and Remark 1]{OttoMenon2004} (see also~\cite{Otto1999}). However, it follows both from the experimental and numerical observations, and pointwise estimates for a reduced model having the same bounds as~\eqref{01-st-transport}-\eqref{01-st-darcy}, that the mixing layer grows as $t$ (see~\cite[Theorem 2 and Theorem 3]{OttoMenon2004} and~\cite[Theorem 1]{OttoMenon2006}). The slow down by a factor of $2$ is termed~\textit{dissipative slowdown}; it is surprising because the speed of traveling waves in viscous scalar conservation laws does not depend on the viscosity. Instead, the  {\em system\/} of PDEs imposes subtle kinematic constraints and appears to have additional dissipative mechanisms that cause the slowdown.

The Rayleigh-Taylor instability is much harder to analyze in similar detail since we lack a comparable understanding of inertial terms and dissipative mechanisms. Our main results capture the optimal scaling of the mixing zone in time. However, a comparison with experimental data reveals that our bounds on the constants $\alpha_\pm$ are an order of magnitude below observations. We expect this gap to be a consequence of the turbulent dissipation. It remains a challenge to obtain sharp bounds on the size of the mixing zone in both the Saffman-Taylor and Rayleigh-Taylor instability.


    \section{The main results}
We state our results in this section.

\subsection{Bulk quantities and notations}
Let us first introduce notation for averaging and normalization. Given a scalar function $f: [0,L]^{d-1} \times \mathbb{R} \to \mathbb{R}$, we average it over the transverse direction
\begin{equation}\label{02-hoiz-average}
    \Bar{f}(z) = \frac{1}{L^{d-1}} \int_{[0,L]^{d-1}} f(y, z) \; dy. 
\end{equation}
In order to ensure $L$-independent estimates, we define normalized integrals
\begin{equation}\label{02-normalized}
    \dashint f(x) \; dx = \int_{\mathbb{R}}\frac{1}{L^{d-1}} \int_{[0,L]^{d-1}} f(y, z) \; dy\,dz.
\end{equation}

The following $L$-independent bulk quantities play a key role in our analysis. First, we consider the negative of the potential energy 
\begin{equation}\label{02-pot-energy}
    E_p(t) = \frac{1}{\rho_+ - \rho_-}\dashint (\rho_0(z) - \rho(x, t)) g z \; dx
\end{equation}
and the kinetic energy of the fluid
\begin{equation}\label{02-kinetic-energy}
    E_k(t) = \frac{1}{2}\frac{1}{\rho_+ - \rho_-} \dashint \rho(x, t) |\u(x, t)|^2 \; dx.
\end{equation}
We also consider the mixing entropies
\begin{equation}
    \label{02-mixing-entropy}
    H(t) = \frac{1}{(\rho_+ - \rho_-)^2} \dashint (\rho_+ - \rho(x, t))(\rho(x, t) - \rho_-) \; dx
\end{equation}
and
\begin{equation}\label{02-mixing-entropy-1}
    \begin{aligned}
        S(t) = & - \dashint \left[ \frac{\rho_+ - \rho(x, t)}{\rho_+ - \rho_-} \log\left( \frac{\rho_+ - \rho(x, t)}{\rho_+ - \rho_-} \right) + \right.& \\ & \left. + \frac{\rho(x, t) - \rho_-}{\rho_+ - \rho_-}  \log\left(\frac{\rho(x, t) - \rho_-}{\rho_+ - \rho_-} \right) \right] \; dx, & 
    \end{aligned}
\end{equation}
which measure the size of ``mushy zones'' where $\rho_- < \rho < \rho_+$. Note that $H(t)$ and $S(t)$ vanish if $\rho \in \{\rho_-, \rho_+\}$ a.e. By the maximum principle, incorporation of the diffusion term in~\eqref{01-transport} ensures $\rho \in (\rho_-,\rho_+)$ for all $x \in [0, L]^{d - 1} \times \mathbb{R}$ and $t > 0$. Thus, $H(t)$ is positive for all $t > 0$. 

Finally, we define the averaged perimeter to measure coarsening
\begin{equation}\label{02-perimeter}
    P(t) = \frac{1}{\rho_+ - \rho_-} \dashint |\nabla \rho(x, t)| \; dx = \frac{1}{\rho_+ - \rho_-}\int_{\rho_-}^{\rho_+} \mathcal{H}^{n - 1} [\rho^{-1} (\eta) ] \; d\eta.
\end{equation}
Here $\mathcal{H}^{n - 1}$ is the $n - 1$ Hausdorff measure and the second equality follows from the co-area formula. Throughout the paper, we denote by $\dot{B}(t)$ the time derivative of a bulk quantity $B(t)$ and by $f_x(x, t)$ and $f_t(x, t)$   the partial derivatives of $f(x, t)$.

We have added a multiplier $(\rho_+ - \rho_-)^{-\alpha}$, $\alpha \in \{1, 2\}$ in the definition of the quantities~\eqref{02-pot-energy}-\eqref{02-perimeter} \textcolor{red} in order to make them invariant under the rescaling of density $\rho'= R \rho$, for any $R>0$. On the other hand, under the spatial rescaling $x'= \lambda x$, we find that the energy, entropies and perimeter scale as $E_p' = \lambda^2 E_p$, $H'=\lambda H$, $S'=\lambda S$ and $P'=P$. The square-root of $E_p/g$ is used to capture the vertical scale of the mixing zone in Theorem~\ref{02-main-thm} below.

\subsection{Uniform bounds on bulk quantities}
\begin{Theorem}\label{02-main-thm}
    Let $(\rho, \u)$ be a classical solution to~\eqref{01-transport}-\eqref{01-momentum} with initial data such that the initial potential energy $E_p(0)$, the mixing entropy $H(0)$, and the mean perimeter $P(0)$ are finite. Then
    \begin{equation}\label{02-energy}
        \limsup_{t \to +\infty} \frac{E_p(t)}{g (A gt^2)^2} \leqslant \frac{1}{24(1 - A^2)}, 
    \end{equation}
    \begin{equation}\label{02-entropy}
        \limsup_{t \to +\infty} \frac{H(t)}{Agt^2} \leqslant \frac{1}{12\sqrt{1 - A^2}},
    \end{equation}
        and
    \begin{equation}
        \limsup_{t \to +\infty} \frac{1}{(Agt^2)^2} \int_0^t P^2(\tau) \; d\tau \leqslant \frac{\pi}{36(1 - A^2)}. \label{02-perimeter-bound}
    \end{equation}
\end{Theorem}

\begin{Remark}
    Theorem~\ref{02-main-thm} provides a comparison between the bounds of the mixing zone and the Riemann problem discussed in Section~\ref{sec:rarefactions}. 
    The bound~\eqref{02-energy} is obtained using the optimal constant in Theorem ~\ref{03-thm-interp}, which is achieved only when the profile corresponds to a rarefaction wave solution to this Riemann problem. This yields the bounds \eqref{eq:scale3} with 
\begin{equation}\label{02-rem-bounds-alpha}
     \alpha_+ =  \frac{1}{1 - A} \quad\mathrm{and}\quad \alpha_- = -\frac{1}{1 + A}.
     \end{equation}
On the other hand, when the mixing zone is modeled as the solution consisting of two shocks to the same Riemann problem we obtain the (formal) bounds
\begin{equation}\label{02-rem-bounds-alpha2}
         \tilde\alpha_+ =  \frac{1}{1 - A + \sqrt{1-A^2}} \quad\mathrm{and}\quad \tilde\alpha_- = -\frac{1}{1 + A +\sqrt{1-A^2}}.
 \end{equation}
     Existence of weak solutions to the {\em immiscible\/} Rayleigh-Taylor instability that mix at this rate is established by Gebhard {\em et al}~\cite{Gebhard2021} provided $A$ is sufficiently close to~$1$. The construction in~\cite{Gebhard2021} also uses a scalar conservation law, but {\em not\/} the one with flux $F(s)$ defined in~\eqref{02-entropy-g} below (see Section ~\ref{sec:rarefactions}).
\end{Remark}


\begin{Remark}
    Relation~\eqref{02-perimeter-bound} is the integrated analog of the (unproven) asymptotic, pointwise inequality
    \begin{equation}\label{02:pertimeter:pointwise}
        P(t) \leqslant \frac{\sqrt{\pi}}{3} \frac{Ag t^{3/2}}{\sqrt{1 - A^2}}.
    \end{equation}
It may be interpreted as a coarsening estimate in the following sense. Define
\begin{equation}
    b_+(t) = \frac{a_+(t)}{P(t)}, \quad b_-(t) = -\frac{a_-(t)}{P(t)},
\end{equation}
and assume that the vertical length scale is much larger than the horizontal length scale, so that
\begin{equation}
    P(t) = \frac{1}{\rho_+ - \rho_-} \dashint |\nabla \rho| \; dx \approx \frac{1}{\rho_+ - \rho_-} \dashint |\partial_y \rho | \; dx.
\end{equation}
Then $b_\pm$ may be considered as an average wavelength of the internal structures. We combine the lower bound on $P(t)$ with the bounds~\eqref{02-bounds-rare} proved below to obtain
\begin{equation}\label{02:horizontal}
    b_+(t) \gtrsim \frac{t^{1/2}}{\sqrt{1 - A}}, \quad b_-(t) \gtrsim \sqrt{1 - A} t^{1/2}.
\end{equation}
\end{Remark}

Let us interpret these estimates physically, though some caution is needed since far less is known about the internal structure of the mixing zone in the Rayleigh-Taylor instability than the Saffman-Taylor instability. The experimental reviews focus on $\alpha_\pm$ and the spreading of the mixing zone; coarsening rates appear not to have been measured systematically (see~\cite{Banerjee2020} for a broad range of experimental diagnostics). This is because it is far harder to probe the internal structure with physical and numerical experiments.

In the Saffman-Taylor instability, the coarsening rates are due to the merging of fingers. By analogy, we expect that in the miscible Rayleigh-Taylor instability, the coarsening is determined by mergers of structures that experimentalists have labeled bubbles and spikes. If the perturbation is more complicated (e.g driven by the white noise), it is not clear what the internal structures are. However, the estimates~\eqref{02:pertimeter:pointwise} and~\eqref{02:horizontal} are robust in the sense that they are free of any assumptions on the form of perturbation. Thus, the interpretation of~\eqref{02:horizontal} is the following: the typical length scale of internal structures must grow at least as $\sqrt{t}$.

\subsection{A crossover effect}\label{subsec:crossover}
Theorem~\ref{02-main-thm} gives an $L$-independent estimate on the potential energy. A strict use of $L$-independent estimates is necessary to obtain self-similar mixing. Our techniques also admit a natural  $L$-dependent modification which describes the late stage of the Rayleigh-Taylor instability. 

\begin{Theorem}\label{02-thm-crossover}
    Let $(\rho, \u)$ be a classical solution to~\eqref{01-transport}-\eqref{01-momentum} with initial data
    such that the initial potential energy $E_p(0)$, the mixing entropy $H(0)$, and the mean
    perimeter $P(0)$ are finite. Then
    \begin{align} \label{02-crossover-estimate}
        \limsup \limits_{t \to +\infty }\frac{E_p(t)}{A^2 g^3 t^3} \lesssim \frac{L^2}{(1 - A)^{3/2}}, \\
        \limsup \limits_{t \to +\infty} \frac{H(t)}{Ag t^{3/2}} \lesssim \frac{L}{(1 - A)^{3/4}}.
    \end{align}
\end{Theorem}

\begin{Remark}
    While the estimates from Theorem~\ref{02-main-thm} are $L$--independent, the estimate~\eqref{02-crossover-estimate} is $L$--dependent because we use a Poincar\'{e} inequality in the transverse direction in the analysis. It provides a different scaling of the size of the mixing layer ($Agt^{3/2}$ instead of $Agt^2$), and becomes sharper when 
    \begin{equation}\label{02-crossover-time}
        t \gg \frac{L^2}{(1 - A)^{1/2}}.
    \end{equation}
    This timescale characterises a regime when dissipation dominates advection and slows down the mixing. 
\end{Remark}

\section{The Riemann problem and some comparisons}\label{sec:rarefactions}

\subsection{Connection with scalar conservation laws}
\label{subsec:Riemann}
Theorem \ref{02-main-thm} allows to compare the classical solution to~\eqref{01-transport}-\eqref{01-momentum} with a solution to a Riemann problem. This idea was introduced in~\cite{Otto1999} and underlies~\cite{OttoMenon2004,OttoMenon2006}. In order to simplify the analysis, we define a function $\chi: [0, L]^{d - 1} \times \mathbb{R} \to [0, 1]$ such that
\begin{equation}\label{02-rho-chi}
    \rho = \rho_- \chi + \rho_+ (1 - \chi),
\end{equation}
and $s: \mathbb{R} \to [0, 1]$ by $s = \Bar{\chi}$. Then
\begin{equation}\label{02-rho-s}
    \Bar{\rho} = \rho_- s + \rho_+ (1 - s).
\end{equation}
The coarse-grained evolution of the ideal mixing zone corresponds to solutions to the scalar conservation law
\begin{equation}\label{02-conservation-law}
    s_\tau + (F(s))_z = 0, \quad \tau = g t^2/2,
\end{equation}
with the concave flux function $F: [0, 1] \to [0, +\infty)$
\begin{equation}\label{02-entropy-g}
    F(s) = (\rho_+ - \rho_-) \frac{s(1 - s)}{\rho_+ s + \rho_-(1 - s)}
\end{equation}
and the Riemann initial data
\begin{equation}\label{02-initial-data}
    s_0 = \begin{cases}
        1, & z \leqslant 0, \\
        0, & z > 0.
    \end{cases}
\end{equation}
This Riemann problem has no classical solutions, but it has infinitely many weak solutions.  In the theory of scalar conservation laws, it is conventional to single out the rarefaction wave as the unique entropy solution to the Riemann problem~\cite{Dafermos}. However, in our work, the Riemann problem arises through an averaging procedure, and this criterion cannot be naively applied. Thus, we distinguish two weak solutions: a rarefaction wave and a two-shock profile. A choice of one of these solutions is motivated by different hypotheses on the structure of the mixing zone.

We now compute these solutions. The rarefaction wave is self-similar and is given in terms of the similarity variable $\xi = z/\tau$ by
\begin{equation}
    s(\xi) = \begin{cases}
        1, &\xi \leqslant -\frac{2A}{1 + A}, \\
        (F')^{-1}(\xi), &-\frac{2A}{1 + A} < \xi \leqslant \frac{2A}{1 - A}, \\
        0, &\xi > \frac{2A}{1 - A},
    \end{cases}
\end{equation}
where  the rarefaction profile is given  by 
\begin{equation}
    (F')^{-1}(\xi) = \frac{1 - A}{2A} \left( \sqrt{\frac{1 + A}{1 - A}}\frac{1}{\sqrt{1 + \xi}} - 1 \right), \; \xi \in \left[-\frac{2A}{1 + A}, \frac{2A}{1 - A} \right].
\end{equation}
This profile gives the following bound of the mixing zone
\begin{align}\label{02-bounds-rare}
    &a_+(t) = \frac{A}{2} \frac{\rho_+ + \rho_-}{\rho_-} gt^2 = \frac{Agt^2}{1 - A}, &a_-(t) = -\frac{A}{2} \frac{\rho_+ + \rho_-}{\rho_+} gt^2 = -\frac{Agt^2}{1 + A},
\end{align}
with the prefactors
\begin{equation}\label{02-bounds-alpha}
    \alpha_+ =  \frac{1}{1 - A} \quad\mathrm{and}\quad \alpha_- = -\frac{1}{1 + A}.
\end{equation}

The flux function $F(s)$ is concave and has its maximum at $s_* = \sqrt{\rho_-}/(\sqrt{\rho_+} + \sqrt{\rho_-})$. Motivated by the proof of Theorem~\ref{03-thm-interp} (see Lemma~\ref{03-lm-monot} for the details), we also define a two-shock solution connecting the states $1$, $s_*$ and $0$ by
\begin{equation}
    s(\xi) = \begin{cases}
        1, &\xi \leqslant -\frac{2A}{1 + A + \sqrt{1 - A^2}}, \\
        \frac{\sqrt{1 - A}}{\sqrt{1 + A} + \sqrt{1 - A}}, &-\frac{2A}{1 + A + \sqrt{1 - A^2}} < \xi \leqslant \frac{2A}{1 - A + \sqrt{1 - A^2}}, \\
        0, &\frac{2A}{1 - A + \sqrt{1 - A^2}} < \xi.
    \end{cases}
\end{equation}
This solution provides the following bounds
\begin{align}
    &\Tilde{a}_+ = \frac{A}{2} \frac{\rho_+ + \rho_-}{\sqrt{\rho_-}(\sqrt{\rho_+} + \sqrt{\rho_-})} gt^2 = \frac{Agt^2}{1 - A + \sqrt{1 - A^2}}, \label{02-bounds-shock-1} \\
    &\Tilde{a}_- = -\frac{A}{2} \frac{\rho_+ + \rho_-}{\sqrt{\rho_+}(\sqrt{\rho_+} + \sqrt{\rho_-})}gt^2 = -\frac{Agt^2}{1 + A + \sqrt{1 - A^2}} \label{02-bounds-shock-2}
 \end{align}
with the prefactors
\begin{equation}\label{02-bounds-alpha2}
    \tilde\alpha_+ =  \frac{1}{1 - A + \sqrt{1-A^2}}, \quad\mathrm{and}\quad \tilde\alpha_- = -\frac{1}{1 + A +\sqrt{1-A^2}}.
\end{equation}

\begin{figure}
    \centering
    \includegraphics[scale=0.25]{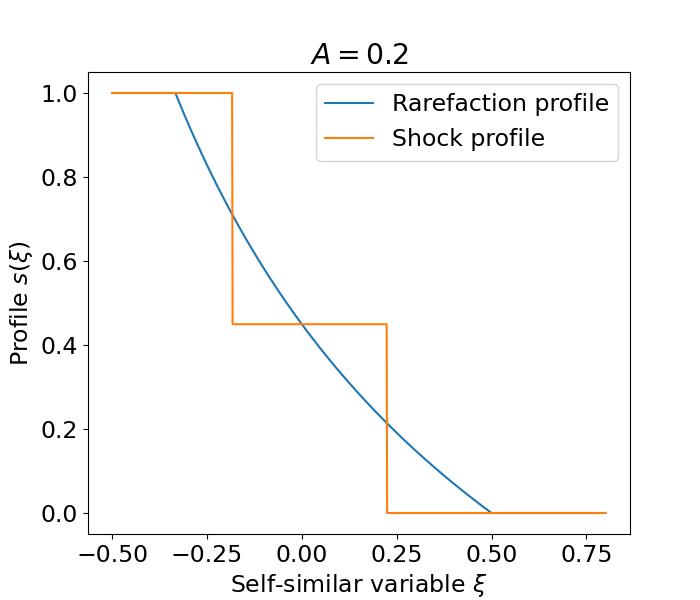}
    \includegraphics[scale=0.25]{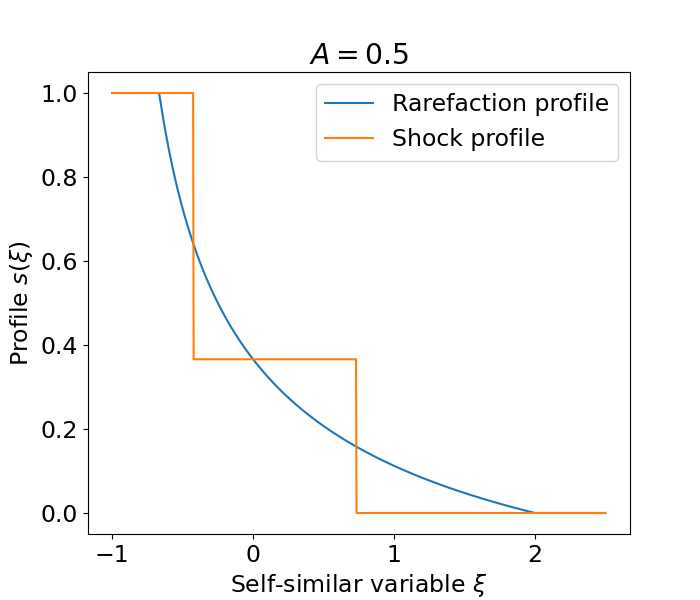}
    \includegraphics[scale=0.25]{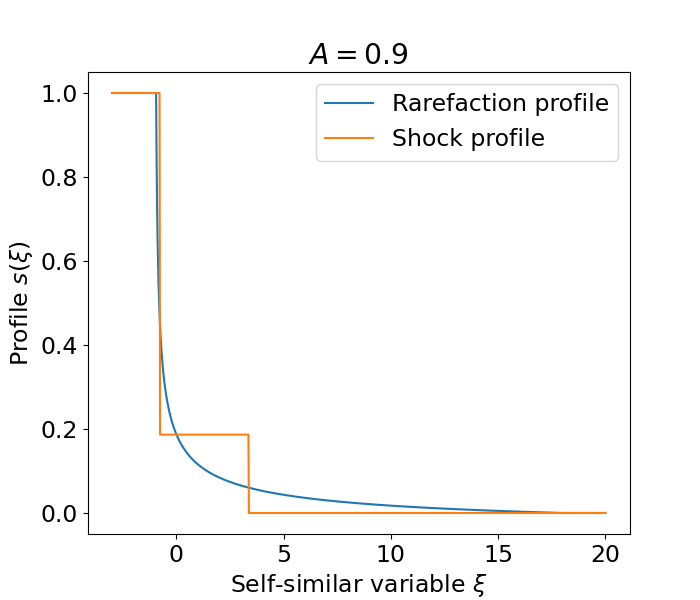}
    \caption{Comparison of the rarefaction (blue) and two-shock (red) profiles in the similarity variable $\xi = x/\tau = 2x/g t^2$ for different values of $A$.}
    \label{02-prifile-comparison}
\end{figure}

The rarefaction and shock profiles are compared in Figure~\ref{02-prifile-comparison}. It immediately follows from the convexity of $F(s)$ that the two-shock solution gives a smaller mixing zone. These two solutions correspond to two different hypotheses on the structure of a homogenized mixing zone. When considering the two-shock solution, we assume that there is a perfectly mixed zone without any $z$ variation. When considering the rarefaction wave we obtain a mixing zone that varies gradually with $z$.

\subsection{Comparison with weak solutions constructed by convex integration}

The $L$-independent upper bounds apply to any classical solution satisfying the hypothesis of Theorem~\ref{02-main-thm}. In particular, any weak solution to the immiscible problem obtained as a limit of a sequence of $L$-dependent classical solutions as $L\to \infty$ must also satisfy these bounds. In contrast, Gebhard~{\em{et al}}~\cite{Gebhard2021} use convex integration to construct weak subsolutions to the {\em immiscible\/} equation
\begin{align}
    &\rho_t + \nabla \cdot (\rho \u) = 0, \label{02-immisc-transport}\\
    &(\rho \u)_t + \nabla \cdot (\rho \u \otimes \u) + \nabla p + \rho g \mathbf{e}_z = 0, &\nabla \cdot \u = 0 \label{02-immisc-momentum},
\end{align}
for the finite time interval $(0, T)$ with $T > 0$ in the finite strip $\Omega = (0, 1) \times (-c_-(T), c_+(T))$.~\footnote{Observe that we study equations~\eqref{01-transport} and~\eqref{01-momentum} on the domain $[0,L]^{d-1}\times \mathbb{R}$. We must rescale space and time, setting $x \mapsto x/L$ and $t \mapsto t/L$, and $g \mapsto Lg$, in order to compare our results with those in~\cite{Gebhard2021}. This rescaling leaves the speed of growth of the mixing zone unchanged, but the diffusion term $\triangle \rho$ is now replaced by $L^{-1}\triangle \rho$ and $-\nabla \rho \cdot \nabla \u$ is now replaced by $-L^{-1} \nabla \rho \cdot \nabla \u$. The rescaled system reduces to equation~\eqref{02-immisc-transport} in the limit $L \to \infty$.} Here the functions $c_\pm(t)$ are given by
\begin{equation}\label{02-c-pm}
    c_\pm(t) = \frac{\rho_+ + \rho_-}{2\sqrt{\rho_\mp}(\sqrt{\rho_+} + \sqrt{\rho_-})} Agt^2,
\end{equation}
and the initial condition is given by perturbation of the flat interface~\eqref{eq:ic1} such that $\rho(x, 0) \in \{ \rho_-, \rho_+\}$ a.e. The main result (see~\cite[Theorems 2.4 and 2.7]{Gebhard2021}) provides the existence of infinitely many bounded admissible weak solutions satisfying $\rho \in \{ \rho_-, \rho_+ \}$ a.e which are mixing within the region $-c_-(t) < z < c_+(t)$, $t \in (0, T)$.

For such weak solutions to the immiscible system~\eqref{02-immisc-transport}-\eqref{02-immisc-momentum}, the transverse averaged profile was obtained in~\cite{Gebhard2021}. If $\Bar{\rho}(z, t)$ is the transverse average and $\tau = g t^2/2$, then $\Bar{\rho}$ satisfies
\begin{equation}
    \Bar{\rho}_{\tau} + (G_0(\Bar{\rho}))_z = 0,
\end{equation}
with the Riemann initial data. The limiting flux function corresponding to the boundary of the convex hull is given by (see \cite[page 1274]{Gebhard2021} for the details)
\begin{equation}
    G_0(\rho) = -\frac{(\rho_+ - \rho)(\rho - \rho_-)}{\rho_+ + \rho_- + \sqrt{\rho_+\rho_-} - \rho}.
\end{equation}
We may also rewrite the flux in terms of $s$ as
\begin{equation}
    G_0(s) =  (\rho_+ - \rho_-)\frac{s(1 - s)}{\rho_+ s + \rho_- (1 - s) + \sqrt{\rho_+ \rho_-}}.
\end{equation}
Taking the entropy solution to the Riemann problem with the flux function $G_0$ provides the bounds given by $c_\pm(t)$ ($\Tilde{a}_\pm$ in our notation). It also gives the following expression for the (non-dimensionalized) negative part of the potential energy~\cite[pp. 1277]{Gebhard2021}
\begin{equation}\label{02-energy-gebhard}
    \begin{split}
        \frac{E_p}{A^2 g^3 t^4} & = \frac{1}{24}\frac{(\rho_+ + \rho_-)^2}{(\sqrt{\rho_+} + \sqrt{\rho_-})^2\sqrt{\rho_+ \rho_-}} \\ & = \frac{1}{24} \frac{2}{\sqrt{1 - A^2}(1 + \sqrt{1 - A^2})} \leqslant \frac{1}{24} \frac{1}{1 - A^2},
    \end{split}
\end{equation}
since $A \in (0, 1)$. The last inequality provides a comparison of bounds from Theorem~\ref{02-main-thm} with the main result from~\cite{Gebhard2021}.

\subsection{Comparison with the experimental data}

We present a comparison of the experimental $\alpha_\pm$ with coefficients obtained from~\eqref{02-bounds-rare} and~\eqref{02-bounds-shock-1}-\eqref{02-bounds-shock-2} in Table~\ref{02-alpha-plus} and Table~\ref{02-alpha-minus}. The experimental and numerical data follows the comprehensive recent review~\cite[Fig.~12]{Banerjee2020}.  We have also compared results from earlier reviews~\cite{Dimonte2000,Dimonte2004} with our work.

\begin{table}[ht]
    \centering
    \begin{tabular}{ || p{2.5cm} | p{2.5cm} | p{2.5cm} | p{2.5cm} ||}\hline
        A & $\alpha_+$ from~\eqref{02-bounds-alpha} & $\alpha_+$ from~\eqref{02-bounds-alpha2} & Observed $\alpha_+$ \\
        \hline
        $0.0025$ & $1.0$ & $0.5$ & $0.065 \pm 0.005$ \\
        $0.08$ & $1.09$ & $0.52$ & $0.065 \pm 0.005$\\
        $0.15$ & $1.18$ & $0.54$ & $0.06$ \\
        $0.2$ & $1.25$ & $0.56$ & $0.055$ \\
        $0.25$ & $1.33$ & $0.58$ & $0.06$ \\
        $0.4$ & $1.67$ & $0.66$ & $0.06$ \\
        $0.46$ & $1.79$ & $0.69$ & $0.06 \pm 0.005$ \\
        $0.6$ & $2.5$ & $0.83$ & $0.06 \pm 0.005$ \\
        $0.72$ & $3.57$ & $1.03$ & $0.04 \pm 0.005$ \\
        $0.8$ & $5.0$ & $1.25$ & $0.05$ \\
        $0.88$ & $8.33$ & $1.68$ & $0.05$ \\
        $0.97$ & $33.33$ & $3.66$ & $0.05$ \\
        \hline
    \end{tabular}
    \caption{Variation of $\alpha_+$ with the Atwood number, $A$.}
    \label{02-alpha-plus}
\end{table}

\begin{table}[ht]
    \centering
    \begin{tabular}{ || p{2.5cm} | p{2.5cm} | p{2.5cm} | p{2.5cm} ||}\hline
        A & $-\alpha_-$ from~\eqref{02-bounds-alpha} & $-\alpha_-$ from~\eqref{02-bounds-alpha2} & Observed $-\alpha_-$ \\
        \hline
        $0.0025$ & $1.0$ & $0.5$ & $0.065 \pm 0.005$  \\
        $0.08$ & $0.93$ & $0.48$ & $0.09 \pm 0.005$ \\
        $0.15$ & $0.87$ & $0.47$ & $0.067 \pm 0.005$ \\
        $0.2$ & $0.83$ & $0.46$ & $0.05 \pm 0.005$ \\
        $0.25$ & $0.8$ & $0.45$ & $0.07 \pm 0.005$ \\
        $0.4$ & $0.71$ & $0.43$ & $0.065$ \\
        $0.46$ & $0.69$ & $0.43$ & $0.08 \pm 0.01$ \\
        $0.6$ & $0.62$ & $0.42$ & $0.09 \pm 0.01$ \\
        $0.72$ & $0.58$ & $0.41$ & $0.075 \pm 0.015$  \\
        $0.8$ & $0.56$ & $0.42$ & $0.11$ \\
        $0.88$ & $0.53$ & $0.42$ & $0.12$ \\
        $0.97$ & $0.51$ & $0.45$ & $0.14$ \\
        \hline
    \end{tabular}
    \caption{Variation of $-\alpha_-$ with the Atwood number, $A$.}
    \label{02-alpha-minus}
\end{table}

The experimental values of $\alpha_\pm$ are significantly less than the predictions. This gap reflects the complexity of the Rayleigh-Taylor instability and may arise from turbulent dissipation, as well as physical effects, including viscosity, that are not modeled by equations~\eqref{01-transport}--\eqref{01-momentum}.

The parameter $\alpha_+$ is consistently smaller than the parameter $-\alpha_-$. This implies that bubbles tend to ascend more rapidly when the prefactor is smaller. Conversely, our estimations suggest $\alpha_+ > -\alpha_-$. This disparity in asymmetry might arise from the influence of various dissipation mechanisms not accounted for in the current system~\eqref{01-transport}-\eqref{01-momentum}. Finally, note that $-\alpha_-(A)$ converges to $g t^2/2$ as $A$ approaches $1$. This limiting scenario corresponds to the free fall of liquid within the gravitational field. 


    \section{Differential inequalities on energy and entropies}\label{sec:diff:ineq}
The first important tool in this paper is differential inequalities for energy and entropies. We start with

\label{sec:proofs}
\begin{Lemma}[Energy balance]\label{03-lm-energy-conservation}
    Let $(\rho, \u)$ be a classical solution to equations~\eqref{01-transport}-\eqref{01-momentum} with the potential energy $E_p(t)$ and the kinetic energy $E_k(t)$ defined in equations~\eqref{02-pot-energy} and~\eqref{02-kinetic-energy} respectively. Then we have
    \begin{align}
        &\dot{E}_p = g - \frac{1}{\rho_+ - \rho_-} \dashint g \rho (\u \cdot \mathbf{e}_z) \; dx, \label{03-e-p-dot} \\ 
        &\dot{E}_k = -\frac{1}{\rho_+ - \rho_-} \dashint g \rho(\u \cdot \mathbf{e}_z) \; dx, \label{03-e-k-dot}
    \end{align}
    and
    \begin{equation}\label{03-energy-conserv}
        E_p(t) - E_k(t) = gt + (E_p(0) - E_k(0)). 
    \end{equation}
\end{Lemma}
\begin{proof}
   We differentiate~\eqref{02-pot-energy} in time and use~\eqref{01-transport} to obtain
    \begin{equation}
        \begin{split}
            \dot{E}_p(t) & = -\frac{1}{\rho_+ - \rho_-}\dashint g z \rho_t \; dx \\ & = - \frac{1}{\rho_+ - \rho_-} \dashint g z (\Delta \rho - \nabla \cdot (\rho \u)) \; dx \\ & = g - \frac{1}{\rho_+ - \rho_-}\dashint g \rho (\u \cdot \mathbf{e}_z) \; dx.
        \end{split}
    \end{equation}
    In the last equality, we use the incompressibility of the velocity field, as well as the periodicity of the boundary conditions concerning $y$. 

    Let us now prove~\eqref{03-e-k-dot}. Taking the time derivative of the kinetic energy, we obtain
    \begin{equation}\label{02-e-k-dot-proof}
        \dot{E}_k = \frac{1}{2} \frac{1}{\rho_+ - \rho_-}\frac{d}{dt}\dashint \rho |\u|^2 \; dx = \frac{1}{\rho_+ - \rho_-}\dashint \left(\rho_t \frac{|\u|^2}{2} + \rho \u \cdot \u_t \right) \; dx
    \end{equation}
    Consider
    \begin{equation}\label{02-energy-flux}
        \begin{split}
            \nabla \cdot \left( \frac{\rho |\u|^2}{2} \u \right) & = \u \cdot \nabla \left( \frac{\rho |\u|^2}{2} \right) \\ & = \u \cdot \left(\rho (\u \cdot \nabla)\u + \frac{|\u|^2}{2} \nabla \rho \right) \\ & = \u \cdot (\rho \u \cdot \nabla)\u + \frac{|\u|^2}{2}(-\rho_t + \Delta \rho).
        \end{split}
    \end{equation}
    In the last equality, we use that $\rho$ satisfies the transport equation~\eqref{01-transport}. Multiply the momentum equation~\eqref{01-momentum} by $\u$ and take the normalized integral,
    \begin{equation}\label{02-a-priori-euler}
        \begin{split}
            \dashint \rho \u \cdot \u_t \; dx & + \dashint \u \cdot (\rho \u \cdot \nabla)\u \; dx \\ & - \dashint (\nabla \rho \cdot \nabla) \u \cdot \u \; dx + \dashint \u \cdot \nabla p \; dx + \dashint g \rho (\u \cdot \mathbf{e}_z) \; dx = 0.
        \end{split}
    \end{equation}
    Note that the third integral may be rewritten as
    \begin{equation}
        -\dashint (\nabla \rho \cdot \nabla) \u \cdot \u \; dx = - \dashint \nabla \rho \cdot \nabla \frac{|\u|^2} {2} \; dx = \dashint \frac{|\u|^2}{2} \Delta \rho \; dx.
    \end{equation}
    The fourth integral in~\eqref{02-a-priori-euler} vanishes due to the incompressibility, and the second integral due to~\eqref{02-energy-flux} is transformed to
    \begin{equation}
        \dashint \u \cdot (\rho \u \cdot \nabla) \u \; dx = -\dashint \frac{|\u|^2}{2}(-\rho_t + \Delta \rho) \; dx.
    \end{equation}
    We finally conclude that
    \begin{equation}\label{02-a-priori-euler-2}
        \dashint \rho \u \cdot \u_t \; dx = -\dashint \frac{|\u|^2}{2} \rho_t\; dx - \dashint g \rho (\u \cdot \mathbf{e}_z) \; dx,
    \end{equation}
    and substituting equation~\eqref{02-a-priori-euler-2} into~\eqref{02-e-k-dot-proof} yields~\eqref{03-e-k-dot}.
\end{proof}

Using the relations above, we can derive the following inequality for the kinetic energy.

\begin{Lemma}\label{03-lm-kinetic}
    Let $(\rho,\u)$ be a classical solution to equations~\eqref{01-transport}-\eqref{01-momentum} with the kinetic energy $E_k(t)$ and $\Tilde{\rho}: \mathbb{R} \to [\rho_-, \rho_+]$ be a background field such that $\rho - \Tilde{\rho} \in L^2$. Then 
    \begin{equation}
        \dot{E}_k \leqslant g\sqrt{2} \left( \frac{1}{\rho_+ - \rho_-} \dashint \frac{(\rho - \Tilde{\rho})^2}{\rho} \; dx \right)^{\frac{1}{2}} E_k^{\frac{1}{2}}.
    \end{equation}
\end{Lemma}
\begin{proof}
    Consider a background field $\Tilde{\rho}: \mathbb{R} \to [\rho_-, \rho_+]$ and note that 
    \begin{equation}\label{03-backgr-field}
        \dashint g \Tilde{\rho}(z) (\u \cdot \mathbf{e}_z) \; dx = 0.
    \end{equation}
    Indeed, taking a horizontal average of the incompressibility equation $\nabla \cdot \u = 0$ and applying periodicity to $\u$ with respect to $y$ yields $(\overline{\u \cdot \mathbf{e}_z})_z = 0$. However, since $|\u| \to 0$ as $z \to \pm \infty$, we have $\overline{(\u \cdot \mathbf{e}_z)} = 0$ for each $t \geqslant 0$. Since $\rho \in [\rho_-, \rho_+]$ with $\rho_- > 0$, $(\rho - \Tilde{\rho})/\sqrt{\rho} \in L^2$. Using~\eqref{03-backgr-field}, one can add $\Tilde{\rho}$ to~\eqref{03-e-k-dot} and then apply the Cauchy-Schwartz inequality,
    \begin{equation}
        \begin{split}
            \dot{E}_k & = - \frac{1}{\rho_+ - \rho_-}\dashint g (\rho - \Tilde{\rho})(\u \cdot \mathbf{e}_z) \; dx \\ & \leqslant \frac{g}{\rho_+ - \rho_-}\left( \dashint \frac{(\rho - \Tilde{\rho})^2}{\rho} \; dx \right)^{\frac{1}{2}} \left( \dashint \rho |\u|^2 \; dx \right)^{\frac{1}{2}}.
        \end{split}
    \end{equation}
    The statement follows directly from the application of the definition of $E_k$.
\end{proof}

We now identify an optimal background field that yields a sharp form of Lemma~\ref{03-lm-kinetic}. Consider
\begin{equation}
    \dashint \frac{(\rho - \Tilde{\rho})^2}{\rho} \; dx = \int_{-\infty}^{+\infty} ( \Bar{\rho} - 2\Tilde{\rho} + \Tilde{\rho}^2 \overline{\rho^{-1}} ) \; dz.
\end{equation}
The minimum is achieved at the point where $-2\Tilde{\rho} + 2\Tilde{\rho} \overline{\rho^{-1}} = 0$. Thus, $\tilde{\rho} = \big(\overline{\rho^{-1}}\big)^{-1}$.

\begin{Lemma}[The flux function]\label{03-lm-optimal-entropy}
    Let $\Tilde{\rho} = \big(\overline{\rho^{-1}} \big)^{-1}$ be the optimal background field. Let $(\rho, \u)$ be a classical solution to equations~\eqref{01-transport}-\eqref{01-momentum}, $\chi: [0, L]^{d - 1} \times \mathbb{R} \to [0, 1]$ and $s: \mathbb{R} \to [0, 1]$ be given by equations~\eqref{02-rho-chi} and~\eqref{02-rho-s}, and $F(s): [0, 1] \to \mathbb{R}$ be given by equation~\eqref{02-entropy-g}. Then
    \begin{equation}\label{03-optimal-entropy-ineq}
        \frac{1}{\rho_+ - \rho_-}\dashint \frac{(\rho - \Tilde{\rho})^2}{\rho} \; dx \leqslant \int_{\mathbb{R}} F(s) \; dz.
    \end{equation}
    The inequality~\eqref{03-optimal-entropy-ineq} is an equality if and only if $\rho \in \{\rho_-, \rho_+\}$ a.e. 
\end{Lemma}
\begin{proof}
    Let $\rho = \rho_- \chi + \rho_+ (1 - \chi)$, and $s = \overline{\chi}$. Applying Jensen's inequality yields
    \begin{equation}
        \frac{1}{\rho_- \chi + \rho_+ (1 - \chi)} \leqslant \frac{\chi}{\rho_-} + \frac{(1 - \chi)}{\rho_+}.
    \end{equation}
    Consequently, 
    \begin{equation}
        \left(\overline{\rho^{-1}}\right)^{-1} \geqslant \frac{\rho_+ \rho_-}{\rho_+ s + \rho_- (1 - s)},
    \end{equation}
    and, hence, 
    \begin{equation}\label{03-average-ineq}
        \overline{\left(\frac{(\rho - \tilde{\rho})^2}{\rho}\right)} = \overline{\rho} - \left(\overline{\rho^{-1}}\right)^{-1} \leqslant (\rho_+ - \rho_-)^2\frac{s(1 - s)}{\rho_+s + \rho_-(1 - s)}.
    \end{equation}
    We divide equation~\eqref{03-average-ineq} by $\rho_+ - \rho_-$ and integrate it over $\mathbb{R}$ to obtain~\eqref{03-optimal-entropy-ineq}. Since $\rho_+ > \rho_- > 0$, $F(s)$ is strictly concave; this, the equality is achieved if and only if $\chi \in \{0, 1\}$ (and $\rho \in \{ \rho_+, \rho_-\}$).
\end{proof}
The next lemma provides the time derivatives for entropies $S(t)$ and $H(t)$.

\begin{Lemma}[Growth of mixing entropies]
    Let $(\rho, \u)$ be a classical solution to~\eqref{01-transport}-\eqref{01-momentum}. Let $H(t)$ and $S(t)$ be the mixing entropies given by~\eqref{02-mixing-entropy} and~\eqref{02-mixing-entropy-1}. Then
    \begin{align}
        &\dot{H} = \frac{2}{(\rho_+ - \rho_-)^2}\dashint |\nabla \rho|^2 \; dx, \label{03-h-dot} \\
        &\dot{S} = \dashint \frac{|\nabla \rho|^2}{(\rho_+ - \rho)(\rho - \rho_-)} \; dx. \label{03-dot-s}
    \end{align}
\end{Lemma}
\begin{proof}
    Let $h: [\rho_-, \rho_+] \to \mathbb{R}$ be a concave function such that $h \in C^2\left((\rho_-, \rho_+)\right)$ and $h(\rho_-) = h(\rho_+) = 0$. Multiply transport equation~\eqref{01-transport} by $h'(\rho)$ and take the normalized integral. Then
    \begin{equation}
        \frac{d}{dt} \dashint h(\rho(x, t)) \; dx + \dashint \u \cdot h'(\rho) \nabla \rho \; dx - \dashint h'(\rho) \Delta \rho \; dx = 0. 
    \end{equation}
    The second term vanishes due to the incompressibility. Indeed,
    \begin{equation}
        \dashint \u \cdot h'(\rho) \nabla \rho \; dx = \dashint \u \cdot \nabla h(\rho) \; dx = \dashint \nabla \cdot (h(\rho) \u) \; dx = 0.
    \end{equation}
   After integrating by parts, the third term is transformed to
    \begin{equation}
        \dashint h'(\rho) \Delta \rho \; dx = -\dashint h''(\rho) |\nabla \rho|^2 \; dx.
    \end{equation}
    Putting
    \begin{align}
        &h_1(\rho) = \frac{(\rho - \rho_-)(\rho_+ - \rho)}{(\rho_+ - \rho_-)^2}, \\
        &h_2(\rho) = - \frac{\rho - \rho_-}{\rho_+ - \rho_-} \log \left( \frac{\rho - \rho_-}{\rho_+ - \rho_-} \right) - \frac{\rho_+ - \rho}{\rho_+ - \rho_-} \log \left( \frac{\rho_+ - \rho}{\rho_+ - \rho_-} \right),
    \end{align}
    finishes the proof.
\end{proof}

\section{A general interpolation inequality}\label{sec:interp:ineq}
In this section, we say $s(z)$ is \textit{a profile} if $s: \mathbb{R} \to [0, 1]$ is measurable satisfying $s(z) \to 1$ as $z \to -\infty$, and $s(z) \to 0$ as $z \to +\infty$. We also denote by $s_0$ the ``stratified'' profile,
\begin{equation}
    s_0(z) = \begin{cases}
        1, & z < 0, \\
        0, & z \geqslant 0.
    \end{cases}
\end{equation}
\begin{Theorem}[Interpolation] \label{03-thm-interp}
    Let $F: [0, 1] \to \mathbb{R}$, $F \in W^{1,2}([0, 1])$ be a concave function such that $F(0) = F(1) = 0$ and let $F$ satisfy the growth condition
    \begin{equation}\label{03-growth-est}
        F(s) \leqslant C_1 s^{\alpha}, \quad F(s) \leqslant C_2(1 - s)^\alpha,
    \end{equation}
    with some constants $C_1, C_2 > 0$ and $\alpha > 1/2$. 
    Then any profile $s$ satisfies the inequality
    \begin{equation}
        \int_{\mathbb{R}} F(s) \; dz \leqslant \mathcal{C} \left( \int_{\mathbb{R}} (s - s_0)z \; dz\right)^{1/2},
    \end{equation}
    with the sharp constant 
\begin{equation}\label{03-c-sharp}
        \mathcal{C}(F) = \left( 2 \int_0^1 (F'(s))^2 \; ds \right)^{1/2}.
    \end{equation}
    Further, the inequality is sharp if and only if
    $s$ is a rarefaction wave profile to the Riemann problem for the hyperbolic conservation law
    \begin{equation}
        s_\tau + (F(s))_z = 0, \quad s(z,0)=s_0(z).
    \end{equation}
\end{Theorem}
\begin{Remark}
    For the flux function $F$ given by~\eqref{02-entropy-g}, the sharp constant is given by
    \begin{equation}\label{03-c-sharp-immisc}
        \mathcal{C}(F) = \sqrt{\frac{2}{3} \frac{(\rho_+ - \rho_-)^2}{\rho_+ \rho_-}}.
    \end{equation}
\end{Remark}

\subsection{Proof of Theorem~\ref{03-thm-interp}}
Proof of the Theorem~\ref{03-thm-interp} is based on several steps. The first step provides the existence of the constant.
\begin{Lemma}[Existence of the constant]\label{03-lm-existence}
    Let $F: [0, 1] \to \mathbb{R}$ be a concave function satisfying conditions from Theorem~\ref{03-thm-interp}. Then there exists a constant $C$ such that for any profile $s$, 
    \begin{equation}
        \int_{\mathbb{R}} F(s) \; dz  \leqslant C \left( \int_{\mathbb{R}} (s - s_0) z \; dz \right)^{1/2}.
    \end{equation}
\end{Lemma}
\begin{proof}
    Fix $R > 0$ and consider the following splitting
    \begin{equation}
        \int_{\mathbb{R}} F(s) \; dz = \int_{-\infty}^{-R} F(s) \; dz + \int_{-R}^R F(s) \; dz + \int_R^{+\infty} F(s) \; dz.
    \end{equation}
    Estimate the first integral. It follows from the growth condition,
    \begin{equation}
        \int_{-\infty}^R F(s) \; dz \leqslant C_1 \int_{-\infty}^R (1 - s)^{\alpha} \; dz = C_1 \int_{-\infty}^{R} (z(s - 1))^{\alpha} (-z)^{-\alpha} \; dz.
    \end{equation}
    Apply the Cauchy-Schwartz inequality and we obtain
    \begin{equation}
        \begin{split}
            \int_{-\infty}^R F(s) \; dz & \leqslant C_1 \int_{-\infty}^R (1 - s)^{\alpha} \; dz \\ & \leqslant C_1\left(\int_{-\infty}^R z(s - 1) \; dz \right)^{\alpha}  \left( \int_{-\infty}^R (-z)^{\alpha/(1 - \alpha)} \; dz \right)^{1 - \alpha} \\ & \leqslant C_1 R^{1 - 2\alpha} \left(\int (s - s_0) z \; dz \right)^{\alpha}.
        \end{split}
    \end{equation}
    Similarly, we obtain the same bound for the third integral. Indeed, if $F(s) < C_2 s^\alpha$ for $\alpha > 1/2$, then
    \begin{equation}
        \int_R^{+\infty} F(s) \; dz \leqslant C_2 R^{1 - 1/2\alpha} \left( \int_\mathbb{R} (s - s_0) z \; dz \right)^{\alpha}.
    \end{equation}
    The integral over the bounded domain is controlled by the $L^\infty$-norm of $F$,
    \begin{equation}
        \int_{-R}^{R} F(s) \; dz \leqslant 2 R \| F \|_{L^\infty},
    \end{equation}
    and taking the optimal value of $R$ provides us
    \begin{equation}
        \int_{\mathbb{R}} F(s) \; dz \leqslant C_3 \| F \|_{L^\infty}^{1 - 1/2\alpha} \left(\int_{\mathbb{R}} (s - s_0) z \; dz \right)^{1/2}.
    \end{equation}
\end{proof}

The next step is to find the optimal value of $\mathcal{C}$. Consider the following variational problem
\begin{equation}\label{03-var-prob}
    \mathcal{C}^2 = \sup_{s} \frac{ \left( \int_{\mathbb{R}} F(s) \; dz \right)^2}{\int_{\mathbb{R}} (s - s_0) z \; dz},
\end{equation}
among all profiles $s$. It follows from Lemma~\ref{03-lm-existence} that the supremum is finite.

First, let us define a rearrangement of a measurable profile $s$. Define functions $s_+(z), s_-(z): \mathbb{R} \to [0, 1]$ by
\begin{equation}
    s_+(z) = \begin{cases}
        s(z), \; &z > 0, \\
        s(-z), \; &z \leqslant 0
    \end{cases}, \quad
    s_-(z) = \begin{cases}
        1 - s(-z), \; &z > 0, \\
        1 - s(z), \; &z \leqslant 0.
    \end{cases}
\end{equation}
Both of the functions $s_\pm(z)$ are measurable and satisfying $s_\pm \to 0$ as $z \to \pm \infty$. Let $s_\pm^*(z)$ be the symmetric decreasing rearrangements of $s_\pm(z)$ (see~\cite{Lieb2001}, Ch.~3 for definition). Define
\begin{equation}
    s_{rearr}(z) = \begin{cases}
        s^*_+(z), \; &z > 0, \\
        1 - s^*_-(z), \; &z \leqslant 0.
    \end{cases}
\end{equation}
Since $s^*_-(z)$ are monotone on $(-\infty, 0)$ and $(0, +\infty)$, $s_{rearr}$ is also monotone on $(-\infty, 0)$ and $(0, +\infty)$. Note that it can still be non-monotone on the whole axis because of the jump at $z = 0$.

\begin{Lemma}[Rearrangement]\label{03-lm-rearr}
    Let $s$ be a profile and $s_{rearr}$ be the corresponding rearrangement defined as above. Then
    \begin{align}
        &\int_\mathbb{R} F(s) \; dz = \int_{\mathbb{R}} F(s_{rearr}(z)) \; dz, \\
        &\int_{\mathbb{R}} (s - s_0) z \; dz \geqslant \int_{\mathbb{R}} (s_{rearr} - s_0) z \; dz.
    \end{align}
\end{Lemma}
\begin{proof}
    Functions $s_\pm$ are symmetric, and the symmetric rearrangement doesn't change the distribution function. It also preserves the distribution function of $s_{rearr}$, and, hence, preserves $\int_{\mathbb{R}} F(s) \; dz$. The second statement follows from the approximation of $s$ by simple functions and applying the rearrangement procedure. 
\end{proof}

It follows from Lemma~\ref{03-lm-rearr}, that rearrangement of a measurable function $s$ increases the constant. Henceforth, we can assume $s(z) = s_{rearr}(z)$.

Let $s$ be a profile monotone on $(-\infty, 0)$ and $(0, +\infty)$. Since $F(s)$ is a strictly convex continuous function, it has a unique maximum in $[0, 1]$. Denote $\zeta = \mathrm{argmax}_{\zeta \in [0, 1]} \; F(\zeta)$ and define
\begin{equation} \label{03-monotonization}
    s_{monot}(z) = \begin{cases}
        \max(s(z), \zeta), &z < 0, \\
        \min(s(z), \zeta), &z > 0,
    \end{cases}
\end{equation}
and we say $s_{monot}(z)$ is a \textit{monotonization} of $s$. Note that $s_{monot}(z)$ is monotone by definition. 

\begin{Lemma}[Monotonization]\label{03-lm-monot}
    Let $s$ be a profile monotone on $(-\infty, 0)$ and $(0, +\infty)$. Let $s_{monot}(z)$  be the corresponding monotonization. Then
    \begin{align}
        &\int_{\mathbb{R}} F(s) \; dz \leqslant \int_{\mathbb{R}} F(s_{monot}(z)) \; dz, \\
        &\int_{\mathbb{R}} (s_{monot} - s_0) z \; dz \leqslant \int_{\mathbb{R}} (s - s_0) z \; dz.
    \end{align}
\end{Lemma}
\begin{proof}
    Since $F(s)$ has a maximum at point $s = \zeta$, for all $z \in \mathbb{R}$,
    \begin{equation}\label{03-pointwise-max}
        F(\min(s(z), \zeta)) \geqslant F(s), \; F(\max(s(z), \zeta)) \geqslant F(s), 
    \end{equation}
    and after the integration we have 
    \begin{equation}
        \begin{split}
            \int_{\mathbb{R}} F(s) \; dz  & = \int_{-\infty}^0 F(s) \; dz + \int_0^{+\infty} F(s) \; dz \\ &\leqslant \int_{-\infty}^0 F(\max(s(z), \zeta)) \; dz + \int_0^{+\infty} F(\min(s(z), \zeta)) \; dz \\ & = \int_{\mathbb{R}} F(s_{monot}(z)) \; dz.
        \end{split}
    \end{equation}
    The second statement can be also obtained from the pointwise definition~\eqref{03-monotonization} after the subtraction of $s_0$, multiplying by $z$ and integration.
\end{proof}

Now we can prove Theorem~\ref{03-thm-interp}.
\begin{proof}
    Due to Lemmas~\ref{03-lm-rearr} and~\ref{03-lm-monot}, the optimal profile $s$ should be monotone. Hence, we can realize it as a probability measure. This approach allows us to find the optimal profile explicitly. It follows from the chain rule,
    \begin{align}
        &\int_\mathbb{R} (s - s_0) z \; dz = \frac{1}{2} \int_0^1 z^2(s) \; ds, \label{03-energy-s}\\
        &\int_\mathbb{R} F(s) \; dz = -\int_0^1 F'(s) z(s) \; ds. \label{03-entropy-s}
    \end{align}
    Applying Cauchy-Schwartz inequality in~\eqref{03-entropy-s} and using~\eqref{03-energy-s} yields
    \begin{equation}
        \begin{split}
            \int F(s) \; dz & = \int_0^1 F'(s) z(s) \; dz \\ & \leqslant \left( 2\int_0^1 (F'(s))^2\; ds \right)^{1/2} \left(\frac{1}{2} \int_0^1 z^2(s) \; ds \right)^{1/2},
        \end{split}
    \end{equation}
    and the inequality becomes an equality if and only if
    \begin{equation}
        \tau F'(s(z)) = z, 
    \end{equation}
    for some $\tau > 0$, i.e $s(z)$ is a rarefaction profile of the scalar conservation law.
\end{proof}

\section{Proofs of the main theorems}
In the section, the potential energy $E_p$ and the kinetic energy $E_k$ are defined by~\eqref{02-pot-energy} and~\eqref{02-kinetic-energy}, the mixing entropies $H(t)$ and $S(t)$ are defined by~\eqref{02-mixing-entropy} and~\eqref{02-mixing-entropy-1}, and the flux function $F(s)$ is defined by~\eqref{02-entropy-g}.
\subsection{Proof of Theorem \ref{02-main-thm}}
It follows from Lemmas~\ref{03-lm-energy-conservation},~\ref{03-lm-kinetic} and~\ref{03-lm-optimal-entropy},
\begin{equation}
    \dot{E}_p \leqslant g + g\sqrt{2}\left(\int_{\mathbb{R}} F(s) \; dz \right)^{1/2} E^{1/2}_k,
\end{equation}
where $F(s)$ is given by~\eqref{02-entropy-g}. Rewrite $E_p$ in terms of $s$,
\begin{equation}
    E_p = \int_{\mathbb{R}} (s - s_0) g z \; dz.
\end{equation}
Apply Theorem \ref{03-thm-interp} and Lemma \ref{03-lm-energy-conservation}, then
\begin{equation}
    \begin{aligned}
        \dot{E}_p & \leqslant g + (2\mathcal{C})^{1/2} g^{3/4} E^{1/4}_p (E_p - gt - E_p(0) + E_k(0))^{1/2} \\ & \leqslant g + (2\mathcal{C})^{1/2} g^{3/4} E_p^{3/4},
    \end{aligned}
\end{equation}
since $t \geqslant 0$, $E_p(t) \geqslant 0$ and $E_k(0) = 0$. Let $e(t): [0, +\infty) \to \mathbb{R}$ solves the following problem
\begin{align}\label{03-energy-bound}
    \dot{e} = g + (2\mathcal{C})^{1/2} g^{3/4} e^{3/4}, \; e(0) = E_p(0).
\end{align}
Then for all $t > 0$, a pointwise inequality $E_p(t) \leqslant e(t)$ holds. The solution to~\eqref{03-energy-bound} is given by
\begin{equation}
    e(t) = \left(\frac{(2\mathcal{C})^{1/2}}{4} g^{3/4} t + E_p(0) \right)^4 + gt.
\end{equation}
Take a limsup of both sides of $E_p(t) \leqslant e(t)$ and plug the explicit expression for $e(t)$, then we get~\eqref{02-energy}.

Now let us prove bounds for $H(t)$ and $S(t)$. One may rewrite them in terms of the function $\chi(x)$,
\begin{align}
    &H(t) = \dashint \chi(1 - \chi) \; dx, \\
    &S(t) = -\dashint (\chi \log \chi + (1 - \chi)\log(1 - \chi)) \; dx.
\end{align}
Let $h: [0, 1] \to \mathbb{R}_+$ be a concave function such that $h(0) = h(1) = 0$. By Jensen's inequality,
\begin{equation}
    \dashint h(\chi(x)) \; dx = \int_{\mathbb{R}} \overline{h(\chi(x))} dx \leqslant \int_{\mathbb{R}} h(s(z)) \; dz.
\end{equation}
Put $h_1(s) = s(1 - s)$ and $h_2(s) = -(s\log s + (1 - s) \log(1 - s))$. Hence, we obtain $H(t) \leqslant C_1 g^{-1/2} E^{1/2}_p$ and $S(t) \leqslant C_2 g^{-1/2} E^{1/2}_p$ with the sharp constants $\mathcal{C}(h_1) = \sqrt{2/3}$ and $\mathcal{C}(h_2) = \pi \sqrt{2/3}$ respectively. substituting the bound for $E_p(t)$ and taking a limsup give bounds for $H(t)$ and $S(t)$.

Now consider the perimeter $P(t)$. Use Cauchy-Schwartz inequality to obtain
\begin{equation}\label{03-p-interp}
    P(t) = \frac{1}{\rho_+ - \rho_-}\dashint |\nabla \rho| \; dx \leqslant \dot{S}^{1/2} H^{1/2}(t).
\end{equation}
Since $S(t)$ is monotone, we get
\begin{equation}\label{03-p-integrated}
    \int_0^t P^2(\tau) \; d\tau = \int_0^t S(\tau) \dot{H}(\tau) \; d\tau \leqslant S(t) H(t).
\end{equation}
Combining bounds for $S(t)$ and $H(t)$ and~\eqref{03-p-integrated} concludes the proof.

\subsection{Proof of Theorem \ref{02-thm-crossover}}
In this case, we do not have to choose the optimal background field. Our goal is to find the timescale when the dissipation mechanism dominates. For this reason, we omit the numerical constants in the inequalities.

Put $\Tilde{\rho} = \Bar{\rho}$ and it follows from Lemma~\ref{03-lm-kinetic},
\begin{equation}\label{03-e-k-crossover}
    \begin{split}
        \dot{E}_k & g\leqslant g \sqrt{\frac{2}{\rho_+ - \rho_-}} \left(\dashint \frac{(\rho - \Bar{\rho})^2}{\rho} \right)^{1/2} E_k^{1/2} \\ 
        & \lesssim \frac{g}{\sqrt{\rho_+ - \rho_-}} \left( \frac{1}{\rho_-} \dashint (\rho - \Bar{\rho})^2 \; dx \right)^{1/2} E_k^{1/2}.
    \end{split}
\end{equation}

For each $z \in \mathbb{R}$ and $t \geqslant 0$, we have $(\rho - \Bar{\rho})(\cdot, z, t) \in L^2([0, L]^{d-1})$. Apply the Poincaré inequality pointwise for all $z$. Since $(\rho - \Bar{\rho})(y, \cdot, t) \in L^2(\mathbb{R})$, we can integrate it over $z$. Hence, we obtain
\begin{equation}
    \dashint (\rho - \Bar{\rho})^2 \; dx \lesssim L^2 \dashint |\rho_y|^2 \; dx \leqslant L^2 \dashint |\nabla \rho|^2 \; dx = L^2 (\rho_+ - \rho_-)^2 \dot{H},
\end{equation}
Plug this inequality to~\eqref{03-e-k-crossover} and we get
\begin{equation}\label{03-e-k-crossover-1}
    \dot{E}_k \lesssim L g \sqrt{\frac{\rho_+}{\rho_-} - 1} \dot{H}^{1/2} E_{k}^{1/2} = Lg \sqrt{\frac{2A}{1 - A}} \dot{H}^{1/2} E_k^{1/2}.
\end{equation}
Divide~\eqref{03-e-k-crossover-1} by $E^{1/2}_k$ and take the time average. Due to the Jensen's inequality,
\begin{equation}
    \begin{split}
        \frac{1}{t} E^{1/2}_k & \lesssim L g \sqrt{\frac{A}{1 - A}} \frac{1}{t} \int_0^t \dot{H}^{1/2} (\tau) \; d\tau  \leqslant \\ & \leqslant L g \sqrt{\frac{A}{1 - A}} \left(\frac{1}{t}\int_0^t \dot{H}(\tau) \; d\tau \right)^{1/2} = L g\sqrt{\frac{A}{1 - A}} \frac{H^{1/2}(t)}{t^{1/2}}.
    \end{split}
\end{equation}
Combining with the estimate for $H(t)$ yields
\begin{equation}
    E_k \lesssim L^2 g^2 \frac{At}{1 - A} H \lesssim L^2 g^2 \frac{At}{1 - A} \frac{A g t^2}{\sqrt{1 - A^2}} \lesssim \frac{L^2 A^2 g^3 t^3}{(1 - A)^{3/2}}.
\end{equation}
Taking $\limsup$ of both sides and applying Lemma~\ref{03-lm-energy-conservation} give the same estimate for $E_p(t)$. The bound on the entropy follows from the the general interpolation inequality.

    \section{Acknowledgements}
      The work of GM was supported by NSF grant DMS-2107205. He is also grateful to the Max Planck Institute for Mathematics in the Sciences at Leipzig and the University of Oslo for hospitality during the completion of this work.

      \section{Conflicts of Interest Statement}
      The authors declare that they have no conflicts of interest related to the work presented in this paper.

    \bibliographystyle{siam}
    \bibliography{sample}

\end{document}